\documentclass[11pt]{article}

\usepackage{epsfig}
\usepackage{latexsym}
\usepackage{amsfonts}

\newtheorem{thm}{Theorem}[section]
\newtheorem{lem}[thm]{Lemma}
\newtheorem{cor}[thm]{Corollary}
\newtheorem{prop}[thm]{Proposition}
\newtheorem{rem}[thm]{Remark}

\newtheorem{conj}[thm]{Conjecture}
\newcommand{\bconj}{\begin{conj}}
\newcommand{\econj}{\end{conj}}
\newcommand{\bth}{\begin{thm}}
\newcommand{\eth}{\end{thm}}
\newcommand{\bl}{\begin{lem}}
\newcommand{\el}{\end{lem}}
\newcommand{\bdf}{\begin{def}}
\newcommand{\edf}{\end{def}}
\newcommand{\bcor}{\begin{cor}}
\newcommand{\ecor}{\end{cor}}
\newcommand{\bprop}{\begin{prop}}
\newcommand{\eprop}{\end{prop}}
\newcommand{\brem}{\begin{rem}}
\newcommand{\erem}{\end{rem}}
\newcommand{\beq}{\begin{equation}}
\newcommand{\eeq}{\end{equation}}
\newcommand{\beqn}{\begin{eqnarray}}
\newcommand{\eeqn}{\end{eqnarray}}
\newcommand{\beqns}{\begin{eqnarray*}}
\newcommand{\eeqns}{\end{eqnarray*}}
\newcommand{\bpr}{\noindent{\bf Proof\hspace{1em}}}
\newcommand{\epr}{\hfill\rule{3mm}{3mm}\vspace{\baselineskip}\\}

\newcommand{\ba}{\begin{array}}
\newcommand{\ea}{\end{array}}

\newcommand{\babs}{\begin{abstract}}
\newcommand{\eabs}{\end{abstract}}

\newcommand{\binom}[2]
{ {#1\choose #2} }

\input{tcilatex}
\input tcilatex

\begin{document}

\title{On the non-3-colourability of random graphs.}
\author{Olivier Dubois \thanks{
LIP6, Box 169, CNRS-Universit\'{e} Paris 6, 4 place Jussieu, 75252
Paris Cedex 05, France.} \and \and
Jacques Mandler\thanks{%
LIP6, Box 169, CNRS-Universit\'{e} Paris 6, 4 place Jussieu, 75252
Paris Cedex 05, France.}}
\date{}
\maketitle


\textbf{Abstract.} We show that for $c\geq 2.4682,$ a random graph
on $n$ vertices with $cn(1+o(1))$ edges almost surely has no
3-colouring. This improves on the current best upper bound of
2.4947.\vspace{0.5in}

\section{Introduction}
An old problem on random graphs remaining open to this day is that
of the existence and determination of a $k$-colourability
threshold, already posed in the paper \cite{ErdRen60} which
launched the whole subject. Using the uniformly distributed model
$G\left( n,m\right) $ of graphs with $m$ edges
on $n$ vertices, it reads: does there exist a constant $c_{k}$ such that if $%
m\sim \left( c_{k}-\varepsilon \right) n$ for some $\varepsilon >0$ as $%
n\rightarrow \infty ,$ then almost all graphs in $G\left(
n,m\right) $ are $k $-colourable, while if $m\sim \left(
c_{k}+\varepsilon \right) n$ almost none is? The conjectured
positive answer has been mostly pursued in the case $k=3,$ the
smallest value for which $k$-colourabilty is an NP-complete
problem. It is now supported by computer experiments which put
$c_{3}$ at about $2.3,$ and a non-uniform version is known to hold
\cite{AchlFrie99}, stating that at least an $n$-dependent
$c_{k}\left( n\right) $ exists with the required property. Whether
$c_{k}\left( n\right) $ converges remains open, but its behaviour
in the large $n$ limit is constrained by proven upper and lower
bounds which are getting progressively tighter.

The first lower bounds were by-products of studies on the existence of a $k$%
-core, which is a necessary condition for non-$k$-colourability, but
understandably perhaps, even the exact threshold for the $k$-core
\cite{PittSpeWor96} yields a mediocre bound for colourability, e.g. $1.675$ for $%
k=3.$ More recently, better bounds were achieved by analyzing
simple colouring algorithms, using the powerful
differential-equation techniques of \cite{Wor95}. This gave
$c_{3}>1.923$ \cite{AchlMol97}, and to go beyond that, more
complex extensions of Wormald's techniques had to be sought
\cite{AchlMoo02}, leading to the current best lower bound of
$2.015$ for $c_{3}.$

Upper bounds are generally based on the ubiquitous first-moment
method, starting with the `naive' bound $c_{k}<k\log k$ as
obtained by Devroye see \cite{Chv91}. For $k=3,$ this is $2.71.$
More recent bounds have used the expected number of colourings
with a local minimality property similar to that introduced for
$k$-SAT \cite{DubBouf96, KirKraKriSta98}. In \cite{DunnZit97},
however,
the property was weaker than it might have been, giving nevertheless $%
c_{3}<2.60.$ This was corrected by \cite{AchlMol99}, introducing
the more restrictive `rigid colourings' but using an untight bound
for a probability appearing in the expectation. This gave
$c_{3}<2.522.$ Finally, this was improved to $2.495$ by
\cite{KapKirSta02}, correcting an error in \cite{KapKirSta00}, and
independently by \cite{FounMcD00}. They used a rather
sophisticated occupancy result of \cite{KamMotPalSpi94} to
evaluate exactly the probability just mentioned.\ We seem to be
reaching a point where the sheer complexity of the calculations
needed to extend the method becomes a hindrance to further
progress.

This paper lowers the upper bound on $c_{3}$ to $2.468155,$ i.e.
by a similar amount as between the current best and previous best.
We still use rigid colourings, but now the graphs themselves are
restricted to a subspace sufficiently `dense' for the first-moment
(Markov) estimate still to apply in the limit, but that also
leaves out many graphs that contribute to raise the first moment.
Perhaps the crucial point is that the calculations are actually
fairly straightforward if well taken. On the other hand, some side
issues may be a bit tedious to check explicitly. In this extended
abstract, we will concentrate on the main flow of the calculation.
We first introduce the subspace of graphs we are restricting our
attention to, and show how calculating the first moment of rigid
colourings in this space yields an upper bound for $c_{3}.\;$We
then perform the actual calculation, leading to a pair of
nonlinear equations in two unknowns. Existence and uniqueness of a
solution is then discussed, justifying a simple iterative
procedure, the result of which is plugged into the expectation.
For the above-quoted value of $c,$ this gives an expectation just
below $1.$

\section{A large subspace of random graphs.}

Throughout the paper, $m=cn\left( 1+o\left( 1\right) \right) ,$ and $%
2.4<c<2.5.$ (Non-colourability is known to hold a.s. above $2.5,$ and what
happens below $2.4$ is not our topic here.) We also denote the probability
space by $G\left( n,m\right) $, and, setting $\lambda =2c,$ by $p\left(
x,\lambda \right) $ or $p_{x}$ the Poisson probability function of mean $%
\lambda ,$ i.e. $e^{-\lambda }\lambda ^{x}/x!$. To avoid irrelevancies in
our enumerations, we will consider directed graphs given by an ordered list
of edges; multiple edges and single-vertex loops are allowed (although
statistically insignificant).\ This changes nothing as to existence and
value of the threshold. With these conventions, $\left| G\left( n,m\right)
\right| =n^{\lambda n}$.

First, let the random variable $\theta _{x}$ denote the number of vertices
of a random graph having degree $x.$ We show that $\theta _{x}$ is
concentrated around its mean which is $p_{x}:$

\begin{lem}
\label{Poiss} There is an absolute constant $C>0$ such that for any $%
\varepsilon >0,$%
\begin{equation}
\Pr \left( \left| \theta _{x}-p_{x}\right| >\varepsilon \right) \leq C\sqrt{%
\lambda n}e^{-c\left( \varepsilon ,p_{x}\right) n},  \label{large_dev1}
\end{equation}
where $c\left( \xi ,\eta \right) =\min \left[ \left( \xi +\eta \right) \log
\left( 1+\xi /\eta \right) -\xi ,\;\xi ^{2}/\left( 2\eta \right) \right] .$
In particular, lim$_{n\rightarrow \infty }\Pr \left( \left| \theta
_{x}-p_{x}\right| >\varepsilon \right) =0.$
\end{lem}

\bpr If $K_{i}$ is the degree of vertex $i,$ the random vector
$\left( K_{i}\right) _{1\leq i\leq n}$ follows a multinomial
distribution, which we can view as describing $\lambda n$
indistinguishable balls (the extremities of the $cn$ edges) being
thrown into $n$ bins (the vertices). A simple Poissonization
argument (considering the situation where a Poisson number $M$ of
balls with mean $\lambda n$ are thrown) shows that there are
independent r.v.'s $L_{i},$ with mean $\lambda ,$ such that the
$K_{i}$'s are
distributed as the $L_{i}$'s, conditional on $M=\lambda n.$ The sum $%
W_{x}^{\prime }=\sum_{1}^{n}\mathbf{1}_{\left\{ L_{i}=x\right\} }$ obeys a
binomial large-deviation inequality:
\begin{eqnarray*}
\Pr \left( \left| \frac{W_{x}^{\prime }}{n}-p_{x}\right| >\varepsilon
\right) \leq 2\;e^{-c\left( \varepsilon ,p_{x}\right) n}
\end{eqnarray*}
where $c\left( \varepsilon ,p_{x}\right) =c_{\varepsilon /p_{x}}p_{x},$ and $%
c_{\varepsilon }$ is as in \cite{AlSp92}, Corollary A.14.
Decomposing w.r.t. the values of $M$ and using a standard
inequality for Poisson r.v.'s yields (\ref {large_dev1}). \epr

It is interesting to note that in the sequel, we do not need the
full strength of Lemma \ref{Poiss}, but only the weak form of
concentration stated at the end.

Now, to $\varepsilon >0$ and $x_{max}\in \Bbb{N,}$ we associate the set $%
\mathcal{G}\left( \varepsilon ,x_{max},n,m\right) $ of graphs such that for $%
0\leq x\leq x_{max},$ the number of vertices with degree $x$ lies between $%
\left( p_{x}-\varepsilon \right) n$ and $\left( p_{x}-\varepsilon \right) n.$
The idea is that in view of Lemma \ref{Poiss}, for $x_{max}$ large enough
and $\varepsilon $ small enough, $\mathcal{G}\left( \varepsilon
,x_{max},n,c\right) $ contains `most well-behaved' members of $G\left(
n,m\right) ,$ and therefore that the expected number of (rigid) $3$%
-colourings of graphs drawn uniformly from $\mathcal{G}\left( \varepsilon
,x_{max},n,m\right) $ provides an upper bound for $\mathbf{Pr} \left( 3\mathrm{-Col}%
\right) ,$ the probability of $3$-colourabiility. But, the
discrepancy between the first-moment (say, rigid-colouring) bound
and $c_{3}$ is due to a 'small' (yet exponential) number of
`rogue' graphs having a huge number of colourings, many of which
are left out from $\mathcal{G}\left( \varepsilon
,x_{max},n,m\right) .$ Consequently, the expectation just
mentioned actually gives a better bound than the rigid colourings
by themselves. Of course, in dealing with $\mathcal{G}\left(
\varepsilon ,x_{max},n,m\right) $ we have to control
approximations.\ As for the calculations, performing them in a
constrained subset of $G\left( n,m\right) $ is actually in some
ways beneficial, e.g. we do not need the balls-and-bins occupancy
results of \cite{KamMotPalSpi94}.

The following proposition makes precise the general idea just explained
while staying at the level of the original probability space $G\left(
n,m\right) ,$ which is technically simpler than working in the subspace $%
\mathcal{G}\left( \varepsilon ,x_{max},n,m\right) .$ The r.v. $R\left(
G\right) $ is defined as the number of rigid colorings of $G\in G\left(
n,m\right) .$

\begin{prop}
\label{epsxM} Let the r.v. $X_{\varepsilon ,x_{max},n,m}$ on $G\left(
n,m\right) $ be defined by
\begin{eqnarray*}
X_{\varepsilon ,x_{max},n,m}\left( G\right) =\left\{
\begin{array}{l}
R\left( G\right) \;\;\;\;\;\mathrm{if}\;\;\;\;\;G\in
\mathcal{G}\left(
\varepsilon ,x_{max},n,m\right)  \\
0\;\;\;\;\;\;\;\;\;\mathrm{otherwise}.
\end{array}
\right.
\end{eqnarray*}
If, for some integer $x_{max}$ and some $\varepsilon >0,\;\;\;\mathbf{E}%
\left[ X_{\varepsilon ,x_{max},n,m}\right] $ tends to $0$ as
$n\rightarrow \infty ,$ then so does $\Pr \left(
3\mathrm{-Col}\right) $.
\end{prop}

\bpr
Let $3$-col$(n,m)=\left\{ G\in G\left( n,m\right) :G\;\;\mathrm{is}\;\;3\mathrm{%
-colourable}\right\} .\;$Then, by Lemma \ref{Poiss} and since $R\left(
G\right) \geq 1$ for $G\in 3$-col$(n,m):$%
\begin{eqnarray*}
\Pr \left( 3\mathrm{-Col}\right)  &=&\frac{\left| 3\mathrm{-col}(n,m)\right| }{%
\left| G\left( n,m\right) \right| } \\
&=&\frac{\left| \mathcal{G}\left( \varepsilon ,x_{max},n,m\right) \cap 3%
\mathrm{-col}(n,m)\right| }{\left| G\left( n,m\right) \right| } \\
&&+\frac{\left| \left\{ G\in G\left( n,m\right) :\exists x\left(
0\leq x\leq x_{max}\;\;\mathrm{and\ \ }\left| \theta _{x}\left(
G\right) -p_{x}\right| \geq \varepsilon \right) \right\} \right|
}{\left| G\left( n,m\right)
\right| } \\
&\leq &\frac{1}{\left| G\left( n,m\right) \right| }\sum_{G\in \mathcal{G}%
\left( \varepsilon ,x_{max},n,m\right) \cap
3\mathrm{-col}(n,m)}1+\left(
x_{max}+1\right) o\left( 1\right)  \\
&\leq &\frac{1}{\left| G\left( n,m\right) \right| }\sum_{G\in \mathcal{G}%
\left( \varepsilon ,x_{max},n,m\right) \cap
3\mathrm{-col}(n,m)}R\left(
G\right) +o\left( 1\right)  \\
&=&\mathbf{E}\left[ X_{\varepsilon ,x_{max},n,m}\right] +o\left( 1\right) .
\end{eqnarray*}
\epr

\section{Combinatorial analysis of the expectation.}

First, let $\Theta _{\varepsilon ,x_{max},n,m}$ be the set of vectors $%
\mathbf{\theta }=\left( \theta _{x}\right) _{0\leq x\leq x_{max}}$ in $%
I_{n}^{x_{max}+1}=\left\{ 0,1/n,2/n,...,1\right\} ^{x_{max}+1}$ with $%
\sum_{x=0}^{x_{max}}\theta _{x}\leq 1,\sum_{x=0}^{x_{max}}x\theta _{x}\leq
\lambda ,$ and all $\left| \theta _{x}-p_{x}\right| <\varepsilon $. Since
this set is of polynomial size $\leq \left( 2\varepsilon n\right)
^{x_{max}+1},$ and polynomial factors are irrelevant in our study, counting
graphs in $\mathcal{G}\left( \varepsilon ,x_{max},n,m\right) $ with some
property $\mathcal{P}$ really boils down to counting, for fixed $\mathbf{%
\theta \in }\Theta _{\varepsilon ,x_{max},n,m},$ graphs with $\mathcal{P}$
in the set $\mathcal{G}\left( \mathbf{\theta }\right) $ of $G\in \mathcal{G}%
\left( \varepsilon ,x_{max},n,m\right) $ with $\left( \theta
_{x}\left( G\right) \right) _{0\leq x\leq x_{max}}=\mathbf{\theta
}.$

To say that $G\in \mathcal{G}\left( \mathbf{\theta }\right) $ is $3$%
-colourable means that there is a partition of the vertices into vertices of
types $0,1$ and $2$ (`blue', `red', `green') such that there are only three
types of edges: types $0$ (joining a blue and a red vertex), $1$ (red and
green), and $2$ (green and blue).

Recall also that a rigid coloring is one in which every vertex of type $0$
has edges joining it to at least a vertex of type $1$ and at least a vertex
of type $2,$ while every vertex of type $1$ is joined to at least a vertex
of type $2.$

Now, given such a partition and rationals $\beta _{0},\beta _{1},\beta
_{2}\in I_{n},$ as well as $\mu _{x,j}^{0},\mu _{x,j}^{1},\mu
_{x,j}^{2}\left( 0\leq x\leq x_{max},0\leq j\leq x\right) ,$ we count the
graphs in $\mathcal{G}\left( \mathbf{\theta }\right) $ with:

\begin{itemize}
\item  $\beta _{i}m$ edges of type $i,$ and:

\item  the vertices of degree $x$ being distributed as follows: for $i=0,1,2$
and $0\leq j\leq x,$ there are $\mu _{x,j}^{i}\theta _{x}n$
vertices of degree $x$ and type $i$, each of which has $j$
(type-$i$) edges joining it to vertices of type $i+1\left(
\mathrm{mod}3\right) ,$
\end{itemize}

which are rigidly coloured. Let $Z\left( \mathbf{\theta },\mathbf{\beta ,\mu
},n,m\right) $ be the number of such graphs.

Note that the colouring being rigid says exactly that $\mu _{x,0}^{0}=\mu
_{x,x}^{0}=0$ and $\mu _{x,0}^{1}=0.\;$Also, any vertex of degree $0$ must
be of type $2,$ and any vertex of degree $1$ must be of type $1$ or $2.$
Accordingly, we will be considering only $\mu _{x,j}^{0}$ for $2\leq x\leq
x_{max},1\leq j\leq x-1,$ and $\mu _{x,j}^{1}$ for $1\leq x\leq
x_{max},1\leq j\leq x.$ For $\mu _{x,j}^{2},$ no restriction applies.

We first choose, among $m$ empty templates representing the edges, those
corresponding to each type of edge, and within each edge template, the
colour of each vertex (recall that in our model the edges are directed).
This can be done in $A_{n}\left( \mathbf{\beta },m\right) $ ways, where
\begin{eqnarray*}
A_{n}\left( \mathbf{\beta },m\right) =\frac{m!}{\left( \beta _{0}m\right)
!\left( \beta _{1}m\right) !\left( \beta _{2}m\right) !}2^{m}.
\end{eqnarray*}
Second, we attribute each vertex a type. Within each group of $\mu
_{x,j}^{i}\theta _{x}n$ vertices of degree $x\leq x_{max}$, we comply with
the above-stated requirements. The remaining $\tau n$ vertices$,$ with $\tau
=1-\sum_{x=0}^{x_{max}}\theta _{x},$ will be those of degree $>x_{max}.$ The
number of ways this can be done is:
\begin{eqnarray*}
B_{n}\left( \mathbf{\theta },\mathbf{\mu },n\right)  &=&\frac{n!}{\left(
\theta _{0}n\right) !\left( \theta _{1}n\right) !...\left( \theta
_{x_{max}}n\right) !\left( \tau n\right) !}\times  \\
&&\prod_{x=0}^{x_{max}}\frac{\left( \theta _{x}n\right) !}{\left( \mu
_{x,1}^{0}\theta _{x}n\right) !...\left( \mu _{x,x-1}^{0}\theta _{x}n\right)
!\left( \mu _{x,1}^{1}\theta _{x}n\right) !...\left( \mu _{x,x}^{1}\theta
_{1}n\right) !\left( \mu _{x,0}^{2}\theta _{x}n\right) !...\left( \mu
_{x,x}^{2}\theta _{x}n\right) !}.
\end{eqnarray*}
Finally, we effectively fill the template locations with the vertices.of
various types. Let $M_{n}\left( \mathbf{\theta },\mathbf{\beta ,\mu }%
,n,m\right) $ be the number of possibilities here. To begin with,
consider the vertices of high degree ($>x_{max}$). They are to
occupy $\sigma n$ places, with $\sigma =\left( \lambda
-\sum_{x=0}^{x_{max}}x\theta _{x}\right) \left( 1+o\left( 1\right)
\right) $ and $\sigma \rightarrow 0$ as $x_{max}\rightarrow
+\infty $ and $\varepsilon \rightarrow 0,$ uniformly in
$\mathbf{\theta }$.\ The ways to assign them are certainly
\textit{less} than
\begin{eqnarray*}
\eta \left( \mathbf{\theta },n,m\right) =\binom{\lambda n}{\sigma n}\left(
\tau n\right) ^{\sigma n}.
\end{eqnarray*}
The $\beta _{0}m$ type-$0$ edge templates contain $2\beta _{1}m$ vertices, $%
\beta _{1}m$ of which are already known to be blue, and $\beta _{1}m$ red.
Let us say that among the blue ones, $\sigma _{0}^{0}m$ are of high degree
and already assigned, and $\sigma _{1}^{0}m$ among the red ones. The number
of ways to fill the still-free places in the type-$0$ templates is, then:
\begin{eqnarray*}
\mathcal{M}_{0}=\frac{\left[ \left( \beta _{0}-\sigma _{0}^{0}\right)
m\right] !}{\prod_{x=2}^{x_{max}}\prod_{j=1}^{x-1}j!^{\mu _{x,j}^{0}\theta
_{x}n}}\frac{\left[ \left( \beta _{0}-\sigma _{1}^{0}\right) \right] m!}{%
\prod_{x=2}^{x_{max}}\prod_{j=1}^{x}\left( x-j\right) !^{\mu
_{x,j}^{1}\theta _{x}n}}.
\end{eqnarray*}
We fill similarly the type-$1$ and type-$2$ templates, with
\begin{eqnarray*}
\mathcal{M}_{1}=\frac{\left[ \left( \beta _{1}-\sigma _{1}^{1}\right)
m\right] !}{\prod_{x=1}^{x_{max}}\prod_{j=1}^{x}j!^{\mu _{x,j}^{1}\theta
_{x}n}}\frac{\left[ \left( \beta _{1}-\sigma _{2}^{1}\right) m\right] !}{%
\prod_{x=0}^{x_{max}}\prod_{j=0}^{x}\left( x-j\right) !^{\mu
_{x,j}^{2}\theta _{x}n}}
\end{eqnarray*}
and
\begin{eqnarray*}
\mathcal{M}_{2}=\frac{\left[ \left( \beta _{2}-\sigma _{2}^{2}\right)
m\right] !}{\prod_{x=0}^{x_{max}}\prod_{j=0}^{x}j!^{\mu _{x,j}^{2}\theta
_{x}n}}\frac{\left[ \left( \beta _{2}-\sigma _{0}^{2}\right) m\right] !}{%
\prod_{x=2}^{x_{max}}\prod_{j=1}^{x-1}\left( x-j\right) !^{\mu
_{x,j}^{0}\theta _{x}n}}
\end{eqnarray*}
possibilities, respectively. By construction $\sigma _{0}^{0}+\sigma
_{1}^{0}+\sigma _{1}^{1}+\sigma _{2}^{1}+\sigma _{2}^{2}+\sigma
_{0}^{2}=\sigma .$ All in all,
\begin{eqnarray*}
M_{n}\left( \mathbf{\theta },\mathbf{\beta ,\mu },n,m\right) \leq \mathcal{M}%
_{0}\mathcal{M}_{1}\mathcal{M}_{2}\eta \left( \mathbf{\theta },n,m\right) ,
\end{eqnarray*}
\begin{eqnarray*}
Z\left( \mathbf{\theta },\mathbf{\beta ,\mu },n,m\right)  &=&A_{n}\left(
\mathbf{\beta },m\right) B_{n}\left( \mathbf{\theta },\mathbf{\mu },n\right)
M_{n}\left( \mathbf{\theta },\mathbf{\beta ,\mu },n,m\right)  \\
&\leq &A_{n}\left( \mathbf{\beta },m\right) B_{n}\left( \mathbf{\theta },%
\mathbf{\mu },n\right) \mathcal{M}_{0}\mathcal{M}_{1}\mathcal{M}_{2}\eta
\left( \mathbf{\theta },n,m\right) .
\end{eqnarray*}
and, for the expectation referred to in Proposition \ref{epsxM}:
\begin{equation}
\mathbf{E}\left[ X_{\varepsilon ,x_{max},n,m}\right] \leq \frac{1}{%
n^{\lambda n}}\left( 2\varepsilon n\right) ^{x_{max}+1}\sum_{\mathbf{\beta
,\mu }}A_{n}\left( \mathbf{\beta },m\right) B_{n}\left( \mathbf{\theta },%
\mathbf{\mu },n\right) M_{n}\left( \mathbf{\theta },\mathbf{\beta ,\mu }%
,n,m\right) ,  \label{start}
\end{equation}
where $\mathbf{\beta }$ and $\mathbf{\mu }$ in the sum are constrained by a
relationship which we shall examine later.

\section{Asymptotics.}

Using an inequality version of Stirling's formula and/or upper bounds on
multinomial coefficients derived from it, it is seen that (with the
convention that an empty product is equal to $1$):
\begin{eqnarray*}
A_{n}\left( \mathbf{\beta },m\right) ^{1/n}\leq \left( 1+o\left( 1\right)
\right) \frac{2^{c}}{\left( \beta _{0}^{\beta _{0}}\beta _{1}^{\beta
_{1}}\beta _{2}^{\beta _{2}}\right) ^{c}},
\end{eqnarray*}
\begin{eqnarray*}
B_{n}\left( \mathbf{\theta },\mathbf{\mu },n\right) ^{1/n}\leq \frac{\left(
1+o\left( 1\right) \right) }{\tau ^{\tau }\prod_{x=0}^{x_{max}}\left( \theta
_{x}\prod_{j=1}^{x-1}\mu _{x,j}^{0\;\;\mu _{x,j}^{0}}\prod_{j=1}^{x}\mu
_{x,j}^{1\;\;\mu _{x,j}^{1}}\prod_{j=0}^{x}\mu _{x,j}^{2\;\;\mu
_{x,j}^{2}}\right) ^{\theta _{x}}},
\end{eqnarray*}
\begin{eqnarray*}
\mathcal{M}_{0}^{1/n} &\leq &\left( 1+o\left( 1\right) \right) \frac{\left(
\beta _{0}cn/e\right) ^{\left( 2\beta _{0}-\sigma _{0}^{0}-\sigma
_{1}^{0}\right) c}}{\prod_{x=0}^{x_{max}}\left( \prod_{j=1}^{x-1}j!^{\mu
_{x,j}^{0}}\prod_{j=1}^{x}\left( x-j\right) !^{\mu _{x,j}^{1}}\right)
^{\theta _{x}}}, \\
\mathcal{M}_{1}^{1/n} &\leq &\left( 1+o\left( 1\right) \right) \frac{\left(
\beta _{1}cn/e\right) ^{\left( 2\beta _{1}-\sigma _{1}^{1}-\sigma
_{2}^{1}\right) c}}{\prod_{x=0}^{x_{max}}\left( \prod_{j=1}^{x}j!^{\mu
_{x,j}^{1}}\prod_{j=0}^{x}\left( x-j\right) !^{\mu _{x,j}^{2}}\right)
^{\theta _{x}}}, \\
\mathcal{M}_{2}^{1/n} &\leq &\left( 1+o\left( 1\right) \right) \frac{\left(
\beta _{2}cn/e\right) ^{\left( 2\beta _{0}-\sigma _{2}^{2}-\sigma
_{0}^{2}\right) c}}{\prod_{x=0}^{x_{max}}\left( \prod_{j=0}^{x}j!^{\mu
_{x,j}^{2}}\prod_{j=1}^{x-1}\left( x-j\right) !^{\mu _{x,j}^{0}}\right)
^{\theta _{x}}},
\end{eqnarray*}
so that, using $j!\left( x-j\right) !=x!/\binom{x}{j}:$%
\begin{eqnarray*}
\left( \mathcal{M}_{0}\mathcal{M}_{1}\mathcal{M}_{2}\right) ^{1/n}\leq \zeta
_{0}\left( \varepsilon ,x_{max}\right) \frac{\left( \beta _{0}^{\beta
_{0}}\beta _{1}^{\beta _{1}}\beta _{2}^{\beta _{2}}\right) ^{2c}\left(
cn/e\right) ^{\left( 2-\sigma \right) c}}{\prod_{x=0}^{x_{max}}x!^{\theta
_{x}}}\prod_{x=0}^{x_{max}}\prod_{j=0}^{x}\binom{x}{j}^{\left( \mu
_{x,j}^{0}+\mu _{x,j}^{1}+\mu _{x,j}^{2}\right) \theta _{x}}
\end{eqnarray*}
where $\lim_{\varepsilon \rightarrow 0,x_{max}\rightarrow +\infty }\zeta
_{0}\left( \varepsilon ,x_{max}\right) =1,$ uniformly in $\mathbf{\theta },$
and that
\begin{eqnarray*}
\left( A_{n}\left( \mathbf{\beta },m\right) B_{n}\left( \mathbf{\theta },%
\mathbf{\mu },n\right) \eta \left( \mathbf{\theta },n,m\right) \mathcal{M}%
_{0}\mathcal{M}_{1}\mathcal{M}_{2}\right) ^{1/n} &\leq &\frac{\zeta
_{0}^{\prime }\left( \varepsilon ,x_{max}\right) }{\prod_{x=0}^{x_{max}}%
\left( x!\theta _{x}\right) ^{\theta _{x}}}\times  \\
&&\frac{2^{c}\left( \beta _{0}^{\beta _{0}}\beta _{1}^{\beta _{1}}\beta
_{2}^{\beta _{2}}\right) ^{c}\left( cn/e\right) ^{2c}}{\prod_{x=0}^{x_{max}}%
\prod_{j=0}^{x}\left\{ \left[ \frac{\mu
_{x,j}^{0}}{\binom{x}{j}}\right] ^{\mu _{x,j}^{0}}\left[ \frac{\mu
_{x,j}^{1}}{\binom{x}{j}}\right] ^{\mu _{x,j}^{1}}\left[ \frac{\mu
_{x,j}^{2}}{\binom{x}{j}}\right] ^{\mu _{x,j}^{2}}\right\}
^{\theta _{x}}},
\end{eqnarray*}
where $\zeta _{0}^{\prime }$ has the same property (and by
convention $0^{0}=1$).

Further, we can get $\mathbf{\theta }$-free estimates where $\theta _{x}$ is
replaced throughout by $p_{x},$ at the price of additional factors which all
tend to $1$ as $\varepsilon \rightarrow 0$ and $x_{max}\rightarrow +\infty
.\;$In the end result, the sum in (\ref{start}) being of a polynomial number
of exponentially-behaved terms,
\begin{eqnarray}
\overline{\lim_{n\rightarrow \infty }}\mathbf{E}\left[ X_{\varepsilon
,x_{max},n,m}\right] ^{1/n} &\leq &\frac{\zeta _{1}\left( \varepsilon
,x_{max}\right) }{2^{c}}\frac{\left( \lambda /e\right) ^{\lambda }}{%
\prod_{x=0}^{x_{max}}\left( x!p_{x}\right) ^{p_{x}}}\times   \label{bas_est}
\\
&&\max_{\mathbf{\beta ,\mu }}\frac{\left( \beta _{0}^{\beta _{0}}\beta
_{1}^{\beta _{1}}\beta _{2}^{\beta _{2}}\right) ^{c}}{\prod_{x=0}^{x_{max}}%
\prod_{j=0}^{x}\left\{ \left[ \frac{\mu
_{x,j}^{0}}{\binom{x}{j}}\right] ^{\mu _{x,j}^{0}}\left[ \frac{\mu
_{x,j}^{1}}{\binom{x}{j}}\right] ^{\mu _{x,j}^{1}}\left[ \frac{\mu
_{x,j}^{2}}{\binom{x}{j}}\right] ^{\mu _{x,j}^{2}}\right\}
^{p_{x}}},\nonumber
\end{eqnarray}
where $\lim_{\varepsilon \rightarrow 0,x_{max}\rightarrow +\infty }\zeta
_{1}\left( \varepsilon ,x_{max}\right) =1,$ and the $\max $ is under
constraints to which we now come.

\section{The optimization.}

We set:
\begin{eqnarray*}
\sum_{j=1}^{x-1}\mu _{x,j}^{0}=\alpha _{x}^{0},\;\;\;\;\;\sum_{j=1}^{x}\mu
_{x,j}^{1}=\alpha _{x}^{1},\;\;\;\;\;\sum_{j=0}^{x}\mu _{x,j}^{0}=\alpha
_{x}^{2},
\end{eqnarray*}
and note that for $0\leq x\leq x_{max},$%
\begin{equation}
\alpha _{x}^{0}+\alpha _{x}^{1}+\alpha _{x}^{2}=1,  \label{Constr-x}
\end{equation}
and that $\alpha _{0}^{2}=\mu _{0,0}^{2}=1,\alpha _{0}^{1}=\mu
_{0,0}^{1}=\alpha _{0}^{0}=\mu _{0,0}^{0}=\alpha _{1}^{0}=\mu _{1,0}^{0}=\mu
_{1,1}^{0}=0;$ also, $\mu _{1,1}^{1}=\alpha _{1}^{1}.$

Introduce the reduced blue, red and green spreads $\varphi _{0},\varphi _{1}$
and $\varphi _{2}$ of the coloured graph.\ These are the quotients by $%
\lambda n$ of the numbers of places of the corresponding colours in our
filled graph template (excluding the places occupied by vertices of degree $%
>x_{max}$); namely,
\begin{equation}
\varphi _{i}=\lambda ^{-1}\sum_{x=0}^{x_{max}}xp_{x}\alpha _{x}^{i}=\lambda
^{-1}\sum_{x=0}^{x_{max}}xp_{x}\sum_{j=0}^{x}\mu _{x,j}^{i}.  \label{def-phi}
\end{equation}
Note that $\beta _{0}+\beta _{2}=2\varphi _{0}+\eta _{0}\left(
\varepsilon ,x_{max}\right) $, where $\lim_{\varepsilon
\rightarrow 0,x_{max}\rightarrow +\infty }\eta _{0}\left(
\varepsilon ,x_{max}\right) =0.$ Similarly, $\beta _{0}+\beta
_{1}=2\varphi _{1}+\eta _{1}$, $\beta _{1}+\beta _{2}=2\varphi
_{2}+\eta _{2}$, and $\varphi _{0}+\varphi _{1}+\varphi
_{2}=1+\eta _{3},$ (actually $\eta _{3}$ depends only on
$x_{max},$ see below). Thus, $\beta _{0}=1-2\varphi _{2}+\eta
_{4},$ and similarly for $\beta _{1}$ and $\beta _{2}.$ From this
we see that (\ref{bas_est}) still holds with $1-2\varphi _{i}$
replacing $\beta _{i+1\left( \mathrm{mod}3\right) },$ and the
maximization on $\mathbf{\mu }$ alone (subject to $0\leq \varphi
_{i}\leq 1/2 $), under the penalty of a slightly larger $\zeta
_{1}\left( \varepsilon ,x_{max}\right) $ which still tends to $1$
in the limit of $\varepsilon \rightarrow 0,x_{max}\rightarrow
+\infty .$ Also, the $\max $ may be extended to $\mathbf{\mu }$
being a vector of reals in $\left[ 0,1\right] .$ We are therefore
looking to solve the problem: minimize the function
\begin{eqnarray*}
f\left( \mathbf{\mu }\right)
&=&\sum_{x=0}^{x_{max}}p_{x}\sum_{j=0}^{x}\left[ \mu _{x,j}^{0}\log \frac{%
\mu _{x,j}^{0}}{\binom{x}{j}}+\mu _{x,j}^{1}\log \frac{\mu _{x,j}^{1}}{%
\binom{x}{j}}+\mu _{x,j}^{2}\log \frac{\mu _{x,j}^{2}}{\binom{x}{j}}\right]
\\
&&+c\left( 1-2\varphi _{0}\right) \log \left( 1-2\varphi _{0}\right)
+c\left( 1-2\varphi _{1}\right) \log \left( 1-2\varphi _{1}\right) +c\left(
1-2\varphi _{2}\right) \log \left( 1-2\varphi _{2}\right)
\end{eqnarray*}
(where by convention $0\log 0=0$), subject to the constraints $%
0=C_{x}=\alpha _{x}^{0}+\alpha _{x}^{1}+\alpha _{x}^{2}-1$ for $0\leq x\leq
x_{max},$ $\mathbf{\mu }\geq \mathbf{0},$ $0\leq \varphi _{0},\varphi
_{1},\varphi _{2}\leq 1/2,$ and $\mu _{x,0}^{0}=\mu _{x,x}^{0}=\mu
_{x,0}^{1}=0$ (so that these are not really variables, and we view $\mathbf{%
\mu }$ as a vector in $\Bbb{R}^{3x_{max}\left( x_{max}+1\right) /2}$ and not
$\Bbb{R}^{3\left( x_{max}+1\right) \left( x_{max}+2\right) /2}$). Setting $%
\varphi _{min}=0.26$ and $\varphi _{max}=0.4,$ simple calculations indicate
that within our chosen $c$ domain, the expected number of (unrestricted) $3$%
-colourings such that $\varphi _{0}\leq \varphi _{min}$ or $\varphi _{0}\geq
\varphi _{max}$ tends to zero anyway, and similarly for $\varphi _{1}$ and $%
\varphi _{2}$ This means that we can restrict $\mathbf{\mu }$ to the set $%
\mathcal{U=}\Bbb{R}_{+}^{3x_{max}\left( x_{max}+1\right) /2}\cap \left\{
\mathbf{\mu }:\varphi _{min}<\varphi _{i}<\varphi _{max},i=0,1,2\right\} .$ $%
\mathcal{U}$ is not open, but it can be seen directly that a vector $\mathbf{%
\mu }$ with a null coordinate (recall that we have excluded the coordinates
\textit{required} to be null) cannot be a local minimum.\ So we can replace $%
\Bbb{R}_{+}$ with $\left] 0,+\infty \right[ $ and minimize on the resulting
open set $\mathcal{D}$, where differential techniques can be used.

Since the constraints above are linear, the classical method of
Lagrange multipliers \cite{Lue84} applies without having to check
for some constraint qualification such as linear independence of
gradients (which is true,
though). Associating a Lagrange multiplier $\Lambda _{x}$ to the constraint $%
C_{x}=0,$ a necessary condition for a local minimum is
\begin{eqnarray*}
0=\nabla f+\sum_{x=0}^{x_{max}}\Lambda _{x}\nabla C_{x}.
\end{eqnarray*}
This gives, for $i=0,1,2,$ with $j\notin \left\{ 0,x\right\} $ if $i=0,$ and
$j\neq 0$ if $i=1:$%
\begin{equation}
\mu _{x,j}^{i}=\binom{x}{j}\frac{\left( 1-2\varphi _{i}\right) ^{x}}{\exp
\left[ \Lambda _{x}/p_{x}+1-x\right] }=\binom{x}{j}\frac{\left( 1-2\varphi
_{i}\right) ^{x}}{\mathcal{B}\left( x,\mathbf{\varphi }\right) }
\label{def-mu}
\end{equation}
which we plug back into $C_{x}=0$ to find that the denominator is
\begin{eqnarray*}
\mathcal{B}\left( x,\mathbf{\varphi }\right) =\max \left( 0,2^{x}-2\right)
\left( 1-2\varphi _{0}\right) ^{x}+\left( 2^{x}-1\right) \left( 1-2\varphi
_{1}\right) ^{x}+2^{x}\left( 1-2\varphi _{2}\right) ^{x}.
\end{eqnarray*}
With the $\varphi _{i}$ defined by (\ref{def-phi}), our necessary condition (%
\ref{def-mu}) is a rather hopeless system of $3x_{max}\left(
x_{max}+1\right) /2$ nonlinear equations in as many unknowns.\ However, its
peculiar form means that if we view (\ref{def-mu}, \ref{def-phi}) as a
system of equations in $\mathbf{\mu }$ and $\mathbf{\varphi ,}$ there is an
`easier' way to solve it, namely eliminating the $\mu $'s by plugging the
r.h.s's $\,r_{x,j}^{i}\left( \mathbf{\varphi }\right) $ of (\ref{def-mu})
into (\ref{def-phi}).\ Noting further that the constraints $C_{x}=0$ imply $%
\varphi _{2}=U\left( x_{max}\right) -\varphi _{0}-\varphi _{1}$, with $%
U\left( x_{max}\right) =\lambda ^{-1}\sum_{x=0}^{x_{max}}xp_{x},$ we obtain
a much nicer necessary condition, namely two equations in the two unknowns $%
\varphi _{0}$ and $\varphi _{1}:$%
\begin{equation}
0=\lambda \varphi
_{0}-\sum_{x=0}^{x_{max}}xp_{x}\sum_{j=0}^{x}r_{x,j}^{0}\left( \varphi
_{0},\varphi _{1},U\left( x_{max}\right) -\varphi _{0}-\varphi _{1}\right) ,
\label{Eq0}
\end{equation}
\begin{equation}
0=\lambda \varphi
_{1}-\sum_{x=0}^{x_{max}}xp_{x}\sum_{j=0}^{x}r_{x,j}^{1}\left( \varphi
_{0},\varphi _{1},U\left( x_{max}\right) -\varphi _{0}-\varphi _{1}\right) .
\label{Eq1}
\end{equation}
(In practice, we take $x_{max}$ sufficiently large so that $U\left(
x_{max}\right) $ can be replaced by $1.$) Having solved (\ref{Eq0}, \ref{Eq1}%
), we recover $\mathbf{\mu }$ from (\ref{def-mu}).

Of course, even a system of two nonlinear equations can be
unmanageable, but in this case a change of variables $\varphi
_{0}=y _{0}+y _{1},\varphi _{1}=y _{0}-y _{1}$ turns (\ref{Eq0},
\ref{Eq1}) into
\begin{equation}
\left\{
\begin{array}{l}
K_{0}\left( y _{0},y _{1}\right) =0, \\
K_{1}\left( y _{0},y _{1}\right) =0,
\end{array}
\right.   \label{Syst1}
\end{equation}
where $K_{0}$ and $K_{1}$ are functions that, within our restricted range $%
\varphi _{min}<\varphi _{i}<\varphi _{max},$ are found to be monotone in
each variable separately, with partial derivatives
\begin{eqnarray*}
\frac{\partial K_{0}}{\partial y _{0}}>0,\;\;\;\frac{\partial K_{0}}{%
\partial y _{1}}>0,\;\;\;\frac{\partial K_{1}}{\partial y _{0}}%
>0,\;\;\;\frac{\partial K_{1}}{\partial y _{1}}<0.
\end{eqnarray*}
Since along $K_{i}=0,$ we have $\partial y _{1}/\partial y
_{0}=-\left( \partial K_{i}/\partial y _{0}\right) /\left(
\partial
K_{i}/\partial y _{1}\right) ,$it follows that $y _{1}$ decreases
in $y _{0}$ along $K_{0}=0,$ while it increases along $K_{1}=0$
(see Fig. 1). Further, for
the smallest attainable value of $y _{0}$ in our range, the solution in $%
y _{1}$ of $K_{0}=0$ is larger than that of $K_{1}=0;$ while the
reverse holds for the largest attainable $y _{0}.\;$ Thus, the
existence and
uniqueness of the solution to (\ref{Syst1}), and therefore also to (\ref{Eq0}%
, \ref{Eq1}), simply follow from the intermediate value theorem. Since a
minimum of $f\left( \mathbf{\mu }\right) $ in $\mathcal{U}$ must exist
(corresponding, in the limit, to the maximum term of $\mathbf{E}\left[
X_{\varepsilon ,x_{max},n,m}\right] $), and since there is no local minimum
on the boundary of $\mathcal{U},$ it must be a point of null gradient in $%
\mathcal{D},$ so this is it.

\begin{figure}[h]
\begin{center}
\vspace{2cm}
\includegraphics[width=100mm,height=60mm]{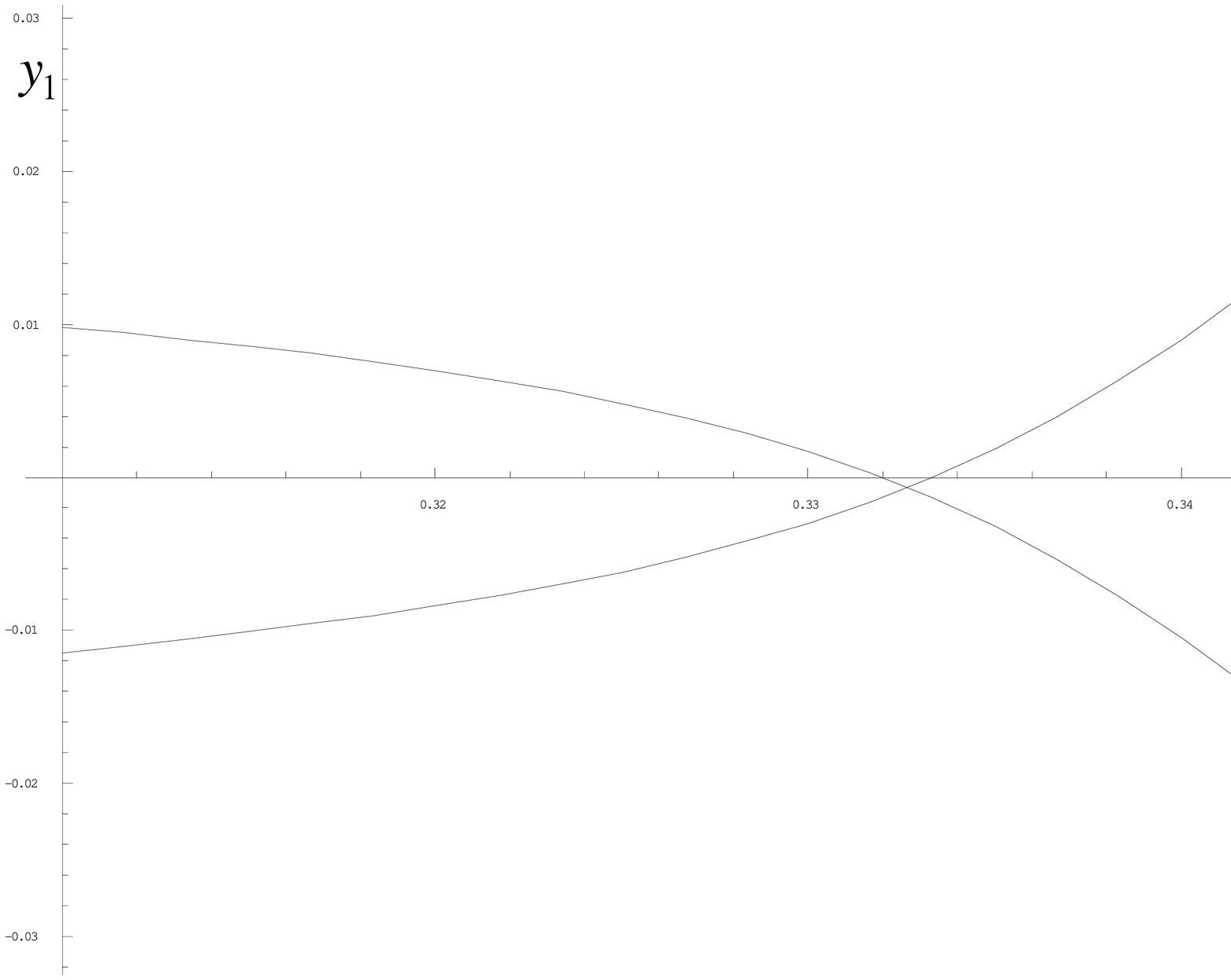}
\vspace{-5cm}
\hspace{+0cm}\parbox{10cm}{\vspace{1.8cm}\caption{The implicit
functions defined by (10). The first equation defines the
decreasing function, the second the increasing one. }}
\label{fig:ratio_heuristique B3}
\end{center}
\end{figure}
\vspace{-5cm}
\newpage
\section{The numerical calculations.}

We now evaluate our modified estimate (\ref{bas_est}) with the
$\varphi _{i}$ and $\mu _{x,j}^{i}$ derived from (\ref{Eq0},
\ref{Eq1}) and (\ref{def-mu}) together with $\mu _{x,0}^{0}=\mu
_{x,x}^{0}=\mu _{x,0}^{1}=0$. We have, for $i=0,1,2$ (since the
factors corresponding to $\mu _{x,j}^{i}=0$ evaluate to $1$ on
both sides of the first equality):
\begin{eqnarray*}
\prod_{x=0}^{x_{max}}\prod_{j=0}^{x}\left[ \frac{\mu _{x,j}^{i}}{\binom{x}{j}%
}\right] ^{\mu _{x,j}^{i}}=\prod_{x=0}^{x_{max}}\prod_{j=0}^{x}\left[ \frac{%
\left( 1-2\varphi _{i}\right) ^{x}}{\mathcal{B}\left( x,\mathbf{\varphi }%
\right) }\right] ^{\mu _{x,j}^{i}}=\prod_{x=0}^{x_{max}}\left[ \frac{\left(
1-2\varphi _{i}\right) ^{x}}{\mathcal{B}\left( x,\mathbf{\varphi }\right) }%
\right] ^{\alpha _{x}^{i}},
\end{eqnarray*}
so, taking account of (\ref{Constr-x}), i.e., $C_{x}=0,$%
\begin{eqnarray*}
\prod_{x=0}^{x_{max}}\prod_{j=0}^{x}\left\{ \left[ \frac{\mu _{x,j}^{0}}{%
\binom{x}{j}}\right] ^{\mu _{x,j}^{0}}\left[ \frac{\mu _{x,j}^{1}}{\binom{x}{%
j}}\right] ^{\mu _{x,j}^{1}}\left[ \frac{\mu _{x,j}^{2}}{\binom{x}{j}}%
\right] ^{\mu _{x,j}^{2}}\right\} ^{p_{x}} &=&\prod_{x=0}^{x_{max}}\frac{%
\left[ \left( 1-2\varphi _{0}\right) ^{\alpha _{x}^{0}}\left( 1-2\varphi
_{1}\right) ^{\alpha _{x}^{1}}\left( 1-2\varphi _{2}\right) ^{\alpha
_{x}^{2}}\right] ^{xp_{x}}}{\mathcal{B}\left( x,\mathbf{\varphi }\right)
^{p_{x}}} \\
&=&\frac{\left( 1-2\varphi _{0}\right) ^{\lambda \varphi _{0}}\left(
1-2\varphi _{1}\right) ^{\lambda \varphi _{1}}\left( 1-2\varphi _{2}\right)
^{\lambda \varphi _{2}}}{\prod_{x=0}^{x_{max}}\mathcal{B}\left( x,\mathbf{%
\varphi }\right) ^{p_{x}}}.
\end{eqnarray*}
Finally, since $\lim_{x_{max}\rightarrow +\infty
}\prod_{x=0}^{x_{max}}\left( x!p_{x}\right) ^{p_{x}}=\left( \lambda
/e\right) ^{\lambda },$%
\begin{eqnarray*}
\overline{\lim_{n\rightarrow \infty }}\mathbf{E}\left[ X_{\varepsilon
,x_{max},n,m}\right] ^{1/n} &\leq &\zeta _{2}\left( \varepsilon
,x_{max}\right) \frac{\prod_{x=0}^{x_{max}}\mathcal{B}\left( x,\mathbf{%
\varphi }\right) ^{p_{x}}}{2^{c}}\times  \\
&&\left( 1-2\varphi _{0}\right) ^{\left( 1-2\varphi _{0}\right) c-\lambda
\varphi _{0}}\left( 1-2\varphi _{1}\right) ^{\left( 1-2\varphi _{1}\right)
c-\lambda \varphi _{1}}\left( 1-2\varphi _{2}\right) ^{\left( 1-2\varphi
_{2}\right) c-\lambda \varphi _{2}} \\
&\leq &\zeta _{2}\left( \varepsilon ,x_{max}\right) \frac{%
\prod_{x=0}^{x_{max}}\mathcal{B}\left( x,\mathbf{\varphi }\right) ^{p_{x}}}{%
2^{c}}\times  \\
&&\left( 1-2\varphi _{0}\right) ^{\left( 1-4\varphi _{0}\right) c}\left(
1-2\varphi _{1}\right) ^{\left( 1-4\varphi _{1}\right) c}\left( 1-2\varphi
_{2}\right) ^{\left( 1-4\varphi _{2}\right) c}
\end{eqnarray*}
where $\lim_{\varepsilon \rightarrow 0,x_{max}\rightarrow +\infty }\zeta
_{2}\left( \varepsilon ,x_{max}\right) =1,$ and an explicit $\zeta
_{2}\left( \varepsilon ,x_{max}\right) $ can be obtained by tracking down
the successive approximations made. Recall that $\varphi _{2}$ is $U\left(
x_{max}\right) -\varphi _{0}-\varphi _{1}$ with $\lim_{x_{max}\rightarrow
+\infty }U\left( x_{max}\right) =1.$

The simple monotonic behaviour described above makes it possible
provably to solve the system (\ref{Eq0}, \ref{Eq1}) using a very
basic iterative procedure that starts from an angle of the
admissible rectangle in $\varphi
_{0},\varphi _{1}$ and spirals towards the solution. Doing so for $%
c=2.468155,$ a sufficiently large $x_{max}$ and small $\varepsilon ,$ one
finds $\overline{\lim_{n\rightarrow \infty }}\mathbf{E}\left[ X_{\varepsilon
,x_{max},n,m}\right] ^{1/n}<0.99999995.$ By monotonicity, this value of $c,$
then, is an upper bound for $c_{3}.\vspace{0.5in}$


\end{document}